\documentclass{amsart}

\usepackage{amsmath}

\usepackage{amssymb}

\usepackage{amsfonts}

\usepackage{amsthm}

\usepackage{latexsym}

\usepackage{url}

\usepackage{verbatim}

\newcommand{\ndash}{\nobreakdash-\hspace{0pt}}

\newcommand{\Ndash}{\nobreakdash--}

\newcommand{\frh}{{\mathfrak{h}}}

\newcommand{\fra}{{\mathfrak{a}}}

\newcommand{\frp}{{\mathfrak{p}}}

\newcommand{\calV}{{\mathcal{V}}}

\newcommand{\calF}{{\mathcal{F}}}

\DeclareMathOperator{\Ker}{Ker}

\DeclareMathOperator{\Der}{Der}

\DeclareMathOperator{\sign}{sign}

\DeclareMathOperator{\dec}{dec}

\newcommand{\id}{\mathrm{id}}

\newcommand{\dd}{\mathrm{d}}

\newtheorem{Theorem}{Theorem}[section]

\newtheorem{Lemma}[Theorem]{Lemma}

\newtheorem*{Theorem*}{Theorem}

\newtheorem{Proposition}[Theorem]{Proposition}

\newtheorem{Corollary}[Theorem]{Corollary}

\theoremstyle{definition}

\newtheorem{Definition}[Theorem]{Definition}

\newtheorem{Example}[Theorem]{Example}

\theoremstyle{remark}

\newtheorem{Remark}[Theorem]{Remark}

\newtheorem*{Ack}{Acknowledgement}

\newcommand{\R}{{\mathbb{R}}}

\newcommand{\N}{{\mathbb{N}}}

\begin{document}

\title{Equivalences of Higher Derived Brackets}

\author{Alberto S. Cattaneo}

\author{Florian Sch\"atz}

\address{Institut f\"ur Mathematik, Universit\"at Z\"urich--Irchel, Winterthurerstrasse 190, CH-8057 Z\"urich, Switzerland}

\email{alberto.cattaneo@math.uzh.ch, florian.schaetz@math.uzh.ch}

%\thanks{The authors acknowledge partial support

%by SNF-grant Nr.20-113439, the European Union through the FP6 Marie

%Curie RTN ENIGMA (contract number MRTN-CT-2004-5652) and by the European Science Foundation through the MISGAM program.}

\thanks{This work has been partially supported by SNF Grant \#20-113439,
by the European Union through the FP6 Marie Curie RTN ENIGMA (contract number MRTN-CT-2004-5652), and by the European Science Foundation through the MISGAM program.}

%\maketitle

\begin{abstract}
This note elaborates on Th.~Voronov's construction \cite{Voronov,Voronov2} of
$L_\infty$\ndash structures via higher derived brackets with a Maurer--Cartan
element. It is shown that gauge equivalent Maurer--Cartan elements induce
$L_\infty$\ndash isomorphic structures. Applications in symplectic, Poisson
and Dirac geometry are discussed.
% This note elaborates on Th.Voronov's constructions of
% $L_\infty$\ndash algebra structures via higher derived brackets constructed with help of an appropriate derivatio given
% in \cite{Voronov} and \cite{Voronov2}.
% It is shown that gauge equivalent derivations induce
% isomorphic $L_{\infty}$\ndash algebra structures. Applications in symplectic, Poisson
% and Dirac geometry are discussed.
\end{abstract}

\maketitle

\section{Introduction}\label{s:intro}

In \cite{Voronov} Th.~Voronov showed that a Maurer--Cartan  element
in a graded Lie algebra  which is split into an abelian subalgebra $\fra$
and another subalgebra $\frp$ induces an $L_\infty$\ndash structure on
the abelian subalgebra $\fra$ in terms of higher derived brackets. 
% In \cite{Voronov} and \cite{Voronov2} Th.Voronov showed that an appropriate derivation
% of a graded Lie algebra $\mathfrak{h}$ which is split into an abelian Lie subalgebra $\mathfrak{a}$
% and another Lie subalgebra $\mathfrak{p}$ induces an $L_\infty$\ndash algebra structure on
% the abelian Lie subalgebra. 

This has interesting applications, e.g., in Poisson 
geometry---especially in view of quantization---where Voronov's construction 
yields an $L_\infty$\ndash structure on the exterior algebra of sections
of the normal bundle of every submanifold (this structure being flat if and only if the submanifold is coisotropic) \cite{OhPark,CattaneoFelder}.
A choice of embedding of the normal bundle is however involved.
It is therefore important
to understand how Voronov's construction depends on this choice. Ultimately this requires understanding how morphisms of graded Lie algebras
influence the induced  $L_\infty$\ndash structures.

It is not difficult to see that
morphisms respecting the splittings 
induce morphisms of the induced $L_\infty$\ndash algebras (see subsection~\ref{s:morph}).
In the application at hand, this implies that a linear automorphism of the normal bundle induces
an $L_\infty$\ndash automorphism (see Remark~\ref{Remark1}).
However, more general diffeomorphisms of the normal bundle do not correspond
to such automorphisms.

The central result of this paper is
that gauge equivalences of Maurer--Cartan elements respecting
the graded Lie subalgebra $\mathfrak{p}$
induce $L_\infty$\ndash automorphisms.
We discuss this $i)$ in the formal setting (Theorem~\ref{formal}) and $ii)$ in case the gauge equivalence is really a flow (Theorem~\ref{MAIN}).
We get an explicit flow, see equations \eqref{U^{1}} and \eqref{U}, of $L_\infty$\ndash algebra automorphisms defined on the
same existence interval.

{}From this we deduce that the $L_\infty$\ndash algebra structure for a submanifold of a Poisson manifold
is canonical up to $L_\infty$\ndash automorphisms (see Section~\ref{s:app}). As a corollary, an isomorphism class of flat
$L_\infty$\ndash algebras  is canonically associated to every regular Dirac manifold (existence of a flat $L_\infty$\ndash structure
was proved in \cite{CZ}).
For the special case of presymplectic manifolds see \cite{OhPark}.

In \cite{Voronov2} it is shown how to extend the original construction to
Maurer--Cartan elements in the graded Lie algebra
of derivations respecting the graded Lie subalgebra $\mathfrak{p}$.
In the present paper we take into account both constructions \cite{Voronov} and \cite{Voronov2}.

\begin{Ack}
We thank J.~Stasheff and M.~Zambon for useful discussions and comments. We also thank the referee
for useful remarks.
\end{Ack}

\section{Higher derived bracket formalism}

We review the higher derived bracket formalism introduced by Th.~Voronov in \cite{Voronov,Voronov2} and explain the problem of finding
`induced automorphisms' in this setting.

\subsection{Preliminaries}\label{s:pre}

Let $V$ be a $\mathbb{Z}$\ndash graded vector space over $\mathbb{R}$ (or any other field of characteristic $0$); i.e.,
$V$ is a collection $\{V_{i}\}_{i \in \mathbb{Z}}$ of vector spaces $V_{i}$ over $\mathbb{R}$.
Homogeneous elements of $V$ of degree $i\in \mathbb{Z}$ are the elements of $V_{i}$. We denote the degree
of a homogeneous element $x \in V$ by $|x|$.  When speaking of linear maps or morphisms, we assume throughout that grading
is preserved. %By morphisms or linear maps between graded vector spaces we mean only degree preserving linear maps.

The $n$th suspension functor $[n]$ from the category of graded vector spaces to itself is defined as follows: given a graded vector space $V$,
$V[n]$ denotes the graded vector space given by the collection $V[n]_{i}:=V_{n+i}$.

One can consider the tensor algebra $T(V)$ associated to a graded vector space $V$ which is a graded vector space
with components
\begin{align*}
T(V)_{m}:=\bigoplus_{k\ge 0}\bigoplus_{j_{1}+\cdots +j_{k}=m}V_{j_{1}}\otimes \cdots \otimes V_{j_{k}}.
\end{align*}
$T(V)$ naturally carries the structure of a cofree coconnected coassociative coalgebra given by the 
deconcatenation coproduct:
\[
\Delta(x_{1}\otimes \dots \otimes x_{n}):=\sum_{i=0}^{n}(x_{1}\otimes \dots \otimes x_{i})\otimes(x_{i+1}\otimes \cdots \otimes x_{n}).
\]
There are two natural representations of the symmetric group $\Sigma_{n}$ on $V^{\otimes n}$: the even one which
is defined by multiplication with the sign $(-1)^{|a||b|}$ for the transposition interchanging $a$ and $b$ in $V$ 
and the odd one by multiplication with the sign $-(-1)^{|a||b|}$ respectively. These two actions naturally 
extend to $T(V)$.
The fix point set of the first action on $T(V)$ is 
denoted by $S(V)$
and called the graded symmetric algebra of $V$ while the fix point set of the latter action is denoted by $\Lambda(V)$
and called the graded skew--symmetric algebra of $V$. The graded symmetric algebra
$S(V)$ inherits a coalgebra structure from $T(V)$ which is cofree
coconnected coassociative and graded cocommutative.
%There are two natural ideals in $T(V)$:
%the one generated by the set $\{x \otimes y - (-1)^{|x||y|}y\otimes x\}$ for all homogeneous $x,y \in V$ which we denote
%by $I_{S}$ and the one
%generated by the set $\{x \otimes y + (-1)^{|x||y|}y\otimes x\}$ which we denote by $I_{\bigwedge}$. We define the
%graded symmetric algebra $S(V)$ to be the quotient of $T(V)$ by $I_{S}$ and the graded skew--symmetric algebra $\bigwedge(V)$
%to be the quotient of $T(V)$ by $I_{\bigwedge}$.

\begin{Definition}
\label{DGLA}
A \textsf{differential graded Lie algebra} $(\mathfrak{h},[\cdot,\cdot])$ is
a graded vector space $\mathfrak{h}$ equipped with a linear map $[\cdot,\cdot]\colon \mathfrak{h}\otimes \mathfrak{h} \to \mathfrak{h}$ 
satisfying the following conditions:
\begin{itemize}
\item graded skew-symmetry: $[x,y]=-(-1)^{|x||y|}[y,x]$,
\item graded Jacobi identity: $[x,[y,z]]=[[x,y],z]+(-1)^{|x||y|}[y,[x,z]]$,
\end{itemize}
for all $x \in \mathfrak{h}_{|x|}$, $y \in \mathfrak{h}_{|y|}$ homogeneous and $z \in \mathfrak{h}$.
\end{Definition}
%\label{gradedLie}
Let $V$ be a graded vector space together with a family of linear maps
\begin{align*}
\{m^{n}\colon S^{n}(V) \to V[1]\}_{n\in \N}.
\end{align*}
Given such a family one defines the associated family of Jacobiators
\begin{align*}
\{J^{n}\colon S^{n}(V) \to V[2]\}_{n \ge 1}
\end{align*}
by
\begin{multline}
\label{Jacobiators}
J^{n}(x_{1} \cdots x_{n}):=\\=
\sum_{r+s=n} \sum_{\sigma \in (r,s)-\text{shuffles}} 
\sign(\sigma)\, m^{s+1}(m^{r}(x_{\sigma(1)} \otimes \cdots \otimes x_{\sigma(r)})
\otimes x_{\sigma(r+1)} \otimes \cdots \otimes x_{\sigma(n)}) 
\end{multline}
where $\sign(\cdot)$ is the Koszul sign, i.e., the one induced from the natural even representation of $\Sigma_{n}$ on $S^{n}(V)$, and
$(r,s)$\ndash shuffles are permutations $\sigma$ of $(1,\dots,r+s)$ such that $\sigma(1) < \cdots < \sigma(r)$ and
$\sigma(r+1) < \cdots < \sigma(n)$.

\begin{Definition}
\label{L-infty}
A family of maps $(m^{n}\colon S^{n}(V) \to V[1])_{n\in \N}$ defines the structure of an \textsf{$L_\infty$\ndash algebra} 
on the graded vector space $V$ whenever the associated family of Jacobiators vanishes identically.
\end{Definition}

This definition is essentially the one given in \cite{Voronov}. We remark that this definition deviates from
the more traditional notion of $L_\infty$\ndash algebras in two points. The early definitions used the graded
skew--symmetric algebra over $V$ instead of the graded symmetric algebra as part of the definition. The transition between
these two settings uses the so called d\'ecalage-isomorphism
\[
\dec^{n}\colon \begin{array}[t]{ccc}
S^{n}(V) &\to& \Lambda^{n}(V[-1])[n]\\
x_{1} \cdots x_{n} &\mapsto& (-1)^{\sum_{i=1}^{n}(n-i)|x_{i}|}x_{1} \wedge \cdots \wedge x_{n}.
\end{array}
\]
%\label{decalage}
More important is the fact that we also allow a `map' $m_{0}\colon \mathbb{R} \to V[1]$ as part of the structure given by an
$L_\infty$\ndash algebra. This piece can be interpreted as an element of $V_{1}$. In the traditional terminology
$m_{0}$ was excluded from the standard definition. Relying on a widespread terminology,
we call structures with $m_{0}=0$ `flat'. Observe that in a flat $L_\infty$\ndash algebra $m_1$ is a differential.

\subsection{V-algebras and induced $L_\infty$-structures}\label{s:ind}

\begin{Definition}
\label{V-algebra}
We call the triple $(\mathfrak{h},\mathfrak{a},\Pi_{\mathfrak{a}})$ a \textsf{V\ndash algebra} (V for Voronov)
if $(\mathfrak{h},[\cdot,\cdot])$ is a graded Lie algebra, $\mathfrak{a}$ is an abelian subalgebra of $\mathfrak{h}$ -- i.e.
$\mathfrak{a}$ is a graded vector subspace of $\mathfrak{h}$ and $[\mathfrak{a},\mathfrak{a}]=0$ -- and
$\Pi_{\mathfrak{a}}\colon \mathfrak{h} \to \mathfrak{a}$ is a projection such that
\begin{align}
\Pi_{\mathfrak{a}}[x,y]=\Pi_{\mathfrak{a}}[\Pi_{\mathfrak{a}}x,y]+\Pi_{\mathfrak{a}}[x,\Pi_{\mathfrak{a}}y]
\label{projection}
\end{align}
holds for every $x, y \in \mathfrak{h}$.
\end{Definition}

Instead of condition \eqref{projection} 
one can require that $\mathfrak{h}$ splits
into $\mathfrak{a}\oplus \mathfrak{p}$ as a vector space where $\mathfrak{p}$ is also a graded Lie subalgebra of $\mathfrak{h}$.
In terms of the projection, $\mathfrak{p}$ is given by the kernel of $\Pi_{\mathfrak{a}}$.

A derivation $E$ of degree $n$ of a graded Lie algebra $\mathfrak{h}$
is a linear map $E \colon \mathfrak{h} \to \mathfrak{h}[n]$ that
satisfies $E[x,y]=[E(x),y]+(-1)^{n|x|}[x,E(y)]$ for all $x \in \mathfrak{h}_{|x|}$, $y\in \mathfrak{h}$.
A derivation $E$ is called inner if there is an element $z \in \mathfrak{h}$ such that $E=[z,\cdot]$.

\begin{Definition}
\label{adaptedder}
Let $(\mathfrak{h},\mathfrak{a},\Pi_{\mathfrak{a}})$ be a V\ndash algebra and $E$ a derivation of $\mathfrak{h}$
that can be written as a sum $E=\Hat E+\Check E$ such that
\begin{itemize}
\item $\Pi_\fra \Hat E\Pi_\fra=\Pi_\fra \Hat E$ (in terms of $\frp:=\Ker\Pi_\fra$ this is equivalent to
$\Hat E(\mathfrak{p})\subset\frp$),
\item $\Check E$ is an inner derivation.
\end{itemize}
Such a derivation $E$ is called \textsf{adapted}. We will denote
the graded Lie algebra of adapted derivations by $\Der(\mathfrak{h},\mathfrak{a},\Pi_{\mathfrak{a}})$.
\end{Definition}

With the help of an adapted derivation $E=\Hat E + [P,\cdot]$ of degree $k$ one can define higher derived brackets on $\mathfrak{a}$:
%Choose a splitting $E=E_{1}+E_{2}$ where $E_{2}:=[P,\cdot]$.
%One defines
\begin{equation}\label{derivedbrackets}
D_{E}^{n}\colon\begin{array}[t]{ccc} \mathfrak{a}^{\otimes n}&\to& \mathfrak{a}[k]\\
x_{1}\otimes \dots \otimes x_{n} &\mapsto& \Pi_{\mathfrak{a}}[[\dots[E(x_{1}),x_{2}],\dots ],x_{n}]
\end{array}
\end{equation}
for every $n > 0$. For $n=0$ we set $D_{E}^{0}:=\Pi_{\mathfrak{a}}P$.
It is easy to check that all these maps are graded commutative; namely,
\begin{multline*}
D_{E}^{n}(x_{1}\otimes \dots \otimes x_{i} \otimes x_{i+1} \otimes \dots \otimes x_{n})=\\
=(-1)^{|x_{i}||x_{i+1}|}D_{E}^{n}(x_{1}\otimes \dots \otimes x_{i+1} \otimes x_{i} \otimes \dots \otimes x_{n})
\end{multline*}
for every $1 \le i \le n-1$.

Observe that for $\Pi_{\mathfrak{p}}:=\id - \Pi_{\mathfrak{a}}$ one can write $E$ as $E=(\Hat E+[\Pi_{\mathfrak{p}}P,\cdot])+[\Pi_{\mathfrak{a}}P,\cdot]$, where
$\Hat E+[\Pi_{\mathfrak{p}}P,\cdot]$ is also a derivation respecting $\mathfrak{p}$, and one 
obtains the same higher derived brackets. So we can always assume without loss of generality that $E$ is the sum
of a derivation respecting $\mathfrak{p}$ and an inner derivation by some element of $\mathfrak{a}$.

We restrict the higher derived brackets constructed from an adapted derivation $E$
to the symmetric algebra $S(\mathfrak{a})$ and obtain a family of maps
%For $E$ an adapted derivation of degree $1$ the higher derived brackets factor through 
%to maps from the symmetric algebra over 
%$\mathfrak{a}$ to $\mathfrak{a}[1]$ and we obtain a family of maps 
$\{D_{E}^{n}\colon S^{n}(\mathfrak{a}) \to \mathfrak{a}[1]\}_{n \in \mathbb{N}}$.

In \cite{Voronov} it is proven that the Jacobiators of the higher derived brackets
$\{D_{E}^{n}\colon S^{n}(\mathfrak{a})\to \mathfrak{a}[1]\}_{n\in \mathbb{N}}$ for $E=[P,\cdot]$ purely inner and of odd degree
are given by the higher derived brackets associated to the inner derivation  associated to $\frac{1}{2}[P,P]$:
\begin{align}
J_{[P,\cdot]}^{n}=D_{[\frac{1}{2}[P,P],\cdot]}^{n}.
\label{J=D}
\end{align}
{}From \eqref{J=D} it follows that all Jacobiators vanish identically
if we assume that $[P,P]=0$ holds. Elements $P$ of degree $1$ that satisfy $[P,P]=0$ are called Maurer--Cartan elements
of $\mathfrak{h}$. Observe that $[\frac{1}{2}[P,P],\cdot]=[P,\cdot]\circ [P,\cdot]$.

In \cite{Voronov2} the case where $E$ is a derivation preserving $\mathfrak{p}$ is considered and it is proved
that for such $E$ of odd degree
\begin{align}
\label{J=D2}
J_{E}^{n}=D_{E\circ E}^{n}
\end{align}
holds. We remark that for an odd derivation $E$, $E\circ E=\frac12[E,E]$ is also a derivation (of even degree). 

This immediately implies that the Jacobiators for any adapted derivation $E$ of odd degree satisfies equation \eqref{J=D2}:
We assume $E=\Hat E+[P,\cdot]$ for $P \in \mathfrak{a}$. One computes
\begin{align*}
J_{E}^{n}=J_{\Hat E}^{n}+D_{[\Hat E(P),\cdot]}^{n}
\end{align*}
and using equation \eqref{J=D2} for $J_{\Hat E}^{n}$ one obtains that equation \eqref{J=D2} holds for all adapted derivations too.
Hence we obtain the following theorem which is a slight variation of similar statements given in \cite{Voronov} and
\cite{Voronov2}:

\begin{Theorem}[Voronov]\label{Voronov}
Let $(\mathfrak{h},\mathfrak{a},\Pi_{\mathfrak{a}})$ be a V\ndash algebra and $E=\Hat E+[P,\cdot]$ 
a Maurer--Cartan element in $\Der(\mathfrak{h},\mathfrak{a},\Pi_{\mathfrak{a}})$.
% be an adapted derivation of $\mathfrak{h}$. 
Then
the family of higher derived brackets associated to $E$,
\begin{align*}
\{D_{E}^{n}\colon S^{n}(\mathfrak{a})\to\mathfrak{a}[1]\}_{n\in \mathbb{N}},
\end{align*}
equips $\mathfrak{a}$ with the structure of an $L_\infty$\ndash algebra in the sense of Definition~\ref{L-infty}. % if $E\circ E=0$ is satisfied.
\end{Theorem}

We remark that the higher derived brackets depend not only on $E$ as a derivation but also on the choice of an element
for the inner derivation. Assume $E=\Hat E+[P,\cdot]$ and $E=\Hat E'+[P',\cdot]$. The two families of
derived brackets for the two decompositions only differ by their $0$\ndash ary operations. In the following we will always assume that
the adapted derivation $E$ comes along with a fixed element $P$ such that $\Hat E+[P,\cdot]$ is the decomposition
of $E$.

\begin{Example}\label{exaDerA}
Let $A$ be a graded commutative algebra and $\Der(A)$ its graded Lie algebra of derivations.
Consider $\frh$ equal to $S_A(\Der(A)[-1])[1]$ or to its formal completion $\Hat S_A(\Der(A)[-1])[1]$.
As $A$ itself is a $\Der(A)$\ndash module, the graded space $\Tilde{\frh}:=A[1]\oplus \Der(A)$ inherits a graded Lie algebra structure.
Since $\Tilde{\frh}[-1]$ generates $\frh[-1]$ as a graded commutative algebra over $A$,
one can extend the Lie bracket uniquely by requiring it to be a graded derivation; namely, one makes $\frh[-1]$ into
a Gerstenhaber algebra. Set $\mathfrak{a}:=A[1]$ and observe that $(\frh,\mathfrak{a},\Pi_{\mathfrak{a}})$ is a V\ndash algebra. Thus, a Maurer--Cartan element induces an $L_\infty$\ndash structure
on $A[1]$ with the additional property that the derived brackets are multiderivations with respect to the 
multiplication in $A$. Such a structure was called $P_\infty$ ($P$ for Poisson) in
\cite{CattaneoFelder}. 

A very special example is when $A=C^\infty(M)$ for a smooth manifold $M$. In this case, $\frh=\calV(M)[1]:=\Gamma(M,\Lambda TM)[1]$ and
the Lie bracket on $\frh$ is the Schouten--Nijenhuis bracket of multivector fields. A Maurer--Cartan element is in this case the same as a Poisson bivector field, and
the induced $P_\infty$\ndash structure is just an ordinary Poisson structure.
More general $P_\infty$\ndash structures are obtained for $M$ a graded manifold.
\end{Example}

\subsection{Morphisms}\label{s:morph}
Suppose now one is given an automorphism $\Phi$ of the graded Lie algebra $\mathfrak{h}$, i.e., a bijective map
$\Phi\colon \mathfrak{h} \to \mathfrak{h}$ that is is degree-preserving and satisfies $\Phi([x,y])=[\Phi(x),\Phi(y)]$ for all
$x, y \in \mathfrak{h}$. If $E$ is a derivation of odd degree, so is $\tilde{E}:=\Phi \circ E \circ \Phi^{-1}$.
Suppose $(\mathfrak{h},\mathfrak{a},\Pi_{\mathfrak{a}})$ is a
V\ndash algebra. One obtains two
families of maps $\{D_{E}^{n}\}_{n \in \mathbb{N}}$ and
$\{D_{\tilde{E}}^{n}\}_{n \in \mathbb{N}}$ that define $L_\infty$\ndash algebra structures on $\mathfrak{a}$. The question
arises under which circumstances these two $L_\infty$\ndash structures are related.

The answer is straightforward as long as the automorphism $\Phi$ respects the splitting.
More generally,
let $\Phi\colon(\frh_1,\fra_1,\Pi_{\fra_1})\to(\frh_2,\fra_2,\Pi_{\fra_2})$ be a morphism of V\ndash algebras, that is,
a morphism of graded Lie algebras $\frh_1\to\frh_2$ satisfying $\Pi_{\fra_2}\circ\Phi=\Phi\circ\Pi_{\fra_1}$. 
Equivalently, $\Phi(\fra_1)\subset\fra_2$ and $\Phi(\frp_1)\subset\frp_2$, with $\frp_i=\Ker\Pi_{\fra_i}$.
We say that $E_i=\Hat E_i + [P_i,\cdot]\in\Der(\frh_i,\fra_i,\Pi_{\fra_i})$, $i=1,2$, are $\Phi$\ndash related
if $E_2\circ\Phi=\Phi\circ E_1$ and $P_2-\Phi(P_1)\in\Ker\Pi_{\fra_2}$.
Then 
\begin{equation}
D_{{E_2}}^{n}(\Phi(x_{1})\otimes\dots \otimes \Phi(x_{n}))=\Phi \circ D_{E_1}^{n}(x_{1}\otimes \dots \otimes x_{n}).
\label{easymorphism}
\end{equation}
Thus, if $E_1$ and $E_2$ are Maurer--Cartan elements, $\Phi$ defines a linear morphism of $L_\infty$\ndash algebras 
$\fra_1\to\fra_2$.

For $E_1=[P_1,\cdot]$ an inner derivation, one may define $E_2= [P_2,\cdot]$ with $P_2=\Phi(P_1)$.
Observe that $E_1$ and $E_2$ are $\Phi$\ndash related and that $E_2$ is Maurer--Cartan if $E_1$ is so.
% \[
% \phi\circ D_{P_1}^k = D_{P_2}^k\circ\phi^{\otimes k},\qquad \forall k\ge0. 
% \]
% As a consequence, if $P_1$ is MC, $\phi$ defines a linear morphism of $L_\infty$\ndash algebras $(\fra_1,D_{P_1})\to(\fra_2,D_{P_2})$.
% This is equivalent to
% $\Pi_{\mathfrak{a}}\circ \Phi=\Phi \circ \Pi_{\mathfrak{a}}$. In this case
% $\Phi$ restricts to an automorphism of $\mathfrak{a}$ and one checks that
% \begin{align}
% D_{\tilde{E}}^{n}(\Phi(x_{1})\otimes\dots \otimes \Phi(x_{n}))=\Phi \circ D_{E}^{n}(x_{1}\otimes \dots \otimes x_{n})
% \label{easymorphism}
% \end{align}
% so $\Phi$ an isomorphism between the two $L_\infty$\ndash algebra structures.

However the requirement on $\Phi$ to respect the splittings
is far too restrictive in general. In the next Section we will show that the conditions under which a family
of automorphisms of $\mathfrak{h}$ induce isomorphisms of the corresponding $L_\infty$\ndash algebras on $\mathfrak{a}$
are much weaker.

\begin{Example}
The V\ndash algebras described in Example~\ref{exaDerA} for $A$ concentrated in degree $0$ (e.g., $A$ the algebra of functions of a smooth manifold)
have the additional property that the splittings respect the degrees (namely, the abelian subalgebra $\mathfrak{a}=A[1]$  and the kernel of the projection are concentrated in degree $-1$
and in nonnegative degrees, respectively). So every graded Lie algebra morphism between such V\ndash algebras is automatically a V\ndash morphism.
\end{Example}

\begin{Example}\label{Example1}
Let $A_1$ and $A_2$ be graded commutative algebras and $\phi\colon A_1\to A_2$ an isomorphism. 
One can extend $\phi$ to an isomorphism of graded Lie algebras $\Phi\colon\Tilde\frh_1:=A_1[1]\oplus \Der(A_1)\to\Tilde\frh_2:=A_2[1]\oplus \Der(A_2)$
by $\Phi(a)=\phi(a)$ for $a\in A_1$ and $\Phi(X)=\phi\circ X\circ \phi^{-1}$ for $X\in\Der(A_1)$. This can be uniquely extended to an isomorphism
$\tilde{\Phi}\colon \frh_1[-1]\to\frh_2[-1]$ of graded commutative algebras, which is also an isomorphism of V\ndash algebras $(\frh_1,A_{1})\to (\frh_2,A_{2})$
(with the canonical projections $\frh_1 \to A_{1}$ and $\frh_2 \to A_{2}$ respectively). 
If we have $\tilde{\Phi}$\ndash related Maurer--Cartan elements, then $\phi$ is an isomorphism of $P_\infty$\ndash algebras.
For example, $\phi$ may be the pushforward of
a diffeomorphism between smooth manifolds or more generally between graded manifolds.
\end{Example}

\section{Induced automorphisms}

Let $(\mathfrak{h},\mathfrak{a},\Pi_{\mathfrak{a}})$ be a V\ndash algebra and
$E=\Hat E + [P,\cdot]$ a Maurer--Cartan element in $\Der(\mathfrak{h},\mathfrak{a},\Pi_{\mathfrak{a}})$.
We denote $\Ker\Pi_\fra$ by $\frp$ throughout. 

The space of Maurer--Cartan elements is invariant under the adjoint action of the
Lie algebra $\Der_0(\mathfrak{h},\mathfrak{a},\Pi_{\mathfrak{a}})$. Such an action is called infinitesimal gauge transformation.
The aim of this Section is to show that integrated gauge transformations preserving $\frp$ induce $L_\infty$\ndash automorphisms.
We do this in the formal and in the analytical setting.
%An infinitesimal gauge equivalence is infinitesimal change of $E$
%by $[m,E]$ with $m\in\Der_0(\mathfrak{h},\mathfrak{a},\Pi_{\mathfrak{a}})$. 
%be an adapted derivation of $\mathfrak{h}$ satisfying
%$E \circ E=0$.

In the formal setting
we introduce a formal parameter $t$ and consider
the V\ndash algebra $(\mathfrak{h}[[t]],\mathfrak{a}[[t]],\Pi_{\mathfrak{a}[[t]]})$ where we use the obvious
$\R[[t]]$\ndash linear extensions of all structure maps. Suppose $m_{t}$ is a derivation of $\mathfrak{h}[[t]]$ of degree $0$.
This derivation
can uniquely be integrated to an automorphism $\phi_{t}$ of $\mathfrak{h}[[t]]$.

In the analytical setting the situation is instead as follows: 
Suppose $m_{t}$ is a family of degree $0$ derivations of $\mathfrak{h}$ for $t \in I$ where $I \subset \mathbb{R}$ is a compact
interval (without loss of generality we will assume that $I=[0,1]$).
We assume that there is a flow $\phi_{t}$ that integrates $m_{t}$ for all $t \in I$.

In both the formal and the analytical setting the flow equation reads
\begin{equation}
\label{Xgeneratingphi}
\begin{aligned}
\frac{\dd}{\dd t}\phi_{t}(z)&= m_{t}\circ\phi_{t}(z),\\
\phi_{0}&=\id,
\end{aligned}
\end{equation}
with the difference that in the formal setting it has to hold for all $z \in \mathfrak{h}[[t]]$ while in the analytical setting
it has to hold for all $z \in \mathfrak{h}$ and all $t\in I$.

We will further assume that
\begin{align}
\label{formalASS_X}
\Pi_{\mathfrak{a}[[t]]}m_{t}\Pi_{\mathfrak{a}[[t]]}=\Pi_{\mathfrak{a}[[t]]}m_t
\end{align}
in the formal setting, and
\begin{align}
\label{ASS_X}
\Pi_{\mathfrak{a}}m_{t}\Pi_{\mathfrak{a}}=\Pi_{\mathfrak{a}}m_{t},\qquad \forall t\in I,
\end{align}
in the analytical setting.

In the formal setting, it follows that the automorphism $\phi_{t}$ satisfies
$\Pi_{\mathfrak{a[[t]]}}\circ \phi_{t} \circ \Pi_{\mathfrak{a[[t]]}}=\Pi_{\mathfrak{a[[t]]}}\circ \phi_{t}$, %respects $\mathfrak{p}[[t]]$, 
while in the analytical setting
the equation 
\begin{align*}
\Pi_{\mathfrak{a}}\circ \phi_{t} \circ \Pi_{\mathfrak{a}}=\Pi_{\mathfrak{a}}\circ \phi_{t}, \qquad \forall t\in I,
\end{align*}
is satisfied under the additional assumption that the only solution to the Cauchy problem
\begin{equation}
\label{ASS_unique}
\begin{aligned}
\frac{\dd}{\dd t}\lambda_{t}&=\Pi_{\mathfrak{a}}m_{t}\lambda_{t},\\
\lambda_{0}&=0,
\end{aligned}
\end{equation}
is $\lambda_{t}=0$ for all $t\in I$.
Equivalently, the condition on $\phi_t$ may be written as
\begin{align}
\label{globalASS_Xformal}
\phi_{t}(\mathfrak{p[[t]]}) &= \mathfrak{p[[t]]}\\
\intertext{and}
\label{globalASS_X}
\phi_{t}(\mathfrak{p}) &= \mathfrak{p}, \qquad \forall t\in I,
\end{align}
respectively.

Finally, in the formal setting, we define $E_{t}:=\phi_{t}\circ \Hat E \circ \phi_{t}^{-1} + [\phi_{t}(P),\cdot]$ and consider
the associated higher derived brackets $\{D_{E_{t}}^{n}\}_{n \in \mathbb{N}}$. Since $\phi_{t}$ satisfies
\eqref{globalASS_Xformal}, $E_{t}$ is an adapted derivation with $E_{t}\circ E_{t}=0$.
Hence we have two $L_\infty$\ndash algebra
structures on $\mathfrak{a}[[t]]$: one is the tautological extension of $\{D_{E}^{n}:S^{n}(\mathfrak{a})\to \mathfrak{a}[1]\}_{n\in \mathbb{N}}$,
which we denote by $\mathfrak{a}[[t]]_{0}$, while the other one is the one associated to $\{D_{E_{t}}^{n}\}_{n \in \mathbb{N}}$, which we denote by $\mathfrak{a}[[t]]_{t}$.
In the analytical setting we consider the one-parameter family of Maurer--Cartan elements
$E_{t}:=\phi_{t}\circ \Hat E \circ \phi_{t}^{-1} + [\phi_{t}(P),\cdot]$
and the associated family of higher derived brackets $\{D_{E_{t}}^{n}\}_{n\in \mathbb{N}}$. We denote
the space $\mathfrak{a}$
equipped with the $L_\infty$\ndash algebra structure defined by the family of maps
$\{D_{E_{t}}^{n}\}_{n \in \mathbb{N}}$ by $\mathfrak{a}_{t}$.

The aim of this Section is to show that, in the formal setting or 
under the condition of uniqueness of solutions to  \eqref{ASS_unique} in the analytical setting, these $L_\infty$\ndash algebra structures
are naturally $L_\infty$\ndash isomorphic. Namely:
%In the analytical setting we additionally will assume the uniqueness for solutions of flow equation \eqref{ASS_unique}.

\begin{Theorem}
\label{formal}
Let $(\mathfrak{h},\mathfrak{a},\Pi_{\mathfrak{a}})$ be a V\ndash algebra and $E$ a Maurer--Cartan element in $\Der(\mathfrak{h},\mathfrak{a},\Pi_{\mathfrak{a}})$.
%Denote the associated $L_\infty$\ndash algebra structure on $\mathfrak{a}[[t]]$ by $\mathfrak{a}[[t]]_{0}$.
Let $\phi_{t}$ be the automorphism of $\mathfrak{h}[[t]]$ generated by a derivation $m_{t}$ of $\mathfrak{h}[[t]]$ of degree $0$
which satisfies \eqref{formalASS_X}.
%Hence $E_{t}:=\phi_{t}\circ E \circ \phi_{t}^{-1}$ is an adapted derivation of $\mathfrak{h}[[t]]$ satisfying $E_{t}\circ E_{t}=0$
%and we denote the
%associated $L_\infty$\ndash algebra structure on $\mathfrak{a}[[t]]$ by $\mathfrak{a}[[t]]_{t}$.
Then the $L_\infty$\ndash algebras $\mathfrak{a}[[t]]_{0}$ and $\mathfrak{a}[[t]]_{t}$ are naturally $L_\infty$\ndash isomorphic.
\end{Theorem}

\begin{Theorem}
\label{MAIN}
Let $(\mathfrak{h},\mathfrak{a},\Pi_{\mathfrak{a}})$ be a V\ndash algebra and $E$ a Maurer--Cartan element in  $\Der(\mathfrak{h},\mathfrak{a},\Pi_{\mathfrak{a}})$.
Assume that $\phi_{t}$ is a family of automorphisms of $\mathfrak{h}$ generated
by a one-parameter family of degree $0$ derivations $m_{t}$ satisfying condition \eqref{ASS_X} and suppose
that equation \eqref{ASS_unique} has a unique solution.
%Consider the family of adapted derivations $E_{t}:=\phi_{t}\circ E \circ \phi_{t}^{-1}$ satisfying $E_{t}\circ E_{t}=0$
%and the
%associated $L_\infty$\ndash algebras $\{\mathfrak{a}_{t}\}_{t \in I}$.
Then the $L_\infty$\ndash algebras $\{\mathfrak{a}_{t}\}_{t\in I}$ are all naturally $L_\infty$\ndash isomorphic.
\end{Theorem}

% We remark that Theorem \ref{MAIN} might eventually hold also in case $\phi_{1}$ is not generated by a flow
% but is an arbitrary automorphism of the graded Lie algebra $\mathfrak{h}$ that respects $\mathfrak{p}$.
% Formula \eqref{U} for the morphism $U$ of $S(\mathfrak{a})$ still
% makes sense in this setting.
The rest of this Section is devoted to the proof of the two Theorems. 
We also get an explicit formula, see \eqref{U^{1}} and \eqref{U}, for the
$L_\infty$\ndash automorphism. Each component of this automorphism is a polynomial in $\phi_t$. So the formula makes sense for every endomorphism
of $\frh$. It is tempting to conjecture that for every graded Lie algebra automorphism respecting $\frp$, it defines an $L_\infty$\ndash automorphism
(or even an $L_\infty$\ndash morphism for every graded Lie algebra endomorphism and a pair of related Maurer--Cartan elements).

\subsection{Infinitesimal considerations}
\label{subsection3.1.}
We briefly review a description of $L_\infty$\ndash algebras, equivalent to the one given in Definition~\ref{L-infty}, which
goes back to Stasheff \cite{Stasheff}. We remarked before that the graded commutative algebra $S(V)$
associated to a graded vector space $V$
is a cofree coconnected graded cocommutative coassociative coalgebra with respect to the coproduct $\Delta$
inherited from $T(V)$.
A linear map $Q\colon S(V) \to S(V)$ that satisfies $\Delta \circ Q= (Q \otimes \text{id} + \text{id} \otimes Q)\circ \Delta$ is called a coderivation
of $S(V)$. By cofreeness of the coproduct $\Delta$ it follows that every linear map from $S(V)$ to $V$
can be extended to a coderivation of $S(V)$ and that every coderivation $Q$ is uniquely determined
by $\text{pr}\circ Q$ where $\text{pr}\colon S(V)\to V$ is the natural projection. So there is a one-to-one
correspondence between families of linear maps $\{m^{n}\colon S^{n}(V) \to V[1]\}_{n\in \mathbb{N}}$ 
and coderivations of $S(V)$ of degree $1$. Moreover, the graded commutator equips
$\text{Hom}(S(V),S(V))$ with the structure of a graded Lie algebra
and this Lie bracket restricts to the subspace of coderivations of $S(V)$. Odd coderivations $Q$ that satisfy $[Q,Q]=0$
are in one-to-one correspondence with families of maps whose associated Jacobiators (see formula \eqref{Jacobiators})
vanish identically. Consequently, Maurer--Cartan elements of the space of coderivations of $S(V)$ correspond
exactly to $L_\infty$\ndash algebra structures on $V$. Since
$Q \circ Q =\frac{1}{2}[Q,Q]=0$, Maurer--Cartan elements of the space of coderivations are exactly the
codifferentials of $S(V)$.

We remark that the approach to $L_\infty$\ndash algebras outlined above makes the notion of $L_\infty$\ndash morphisms
especially transparent: these are just coalgebra morphisms that are
chain maps between the graded symmetric algebras equipped with the codifferentials
that define the $L_\infty$\ndash algebra structures.

In particular, we can interpret the $L_\infty$\ndash algebra structure on $\mathfrak{a}[[t]]$ as a codifferential
$Q(t)$ of $S(\mathfrak{a}[[t]])$. In the analytical setting, we interpret the one-parameter family of
$L_\infty$\ndash algebras $\{a_{t}\}_{t \in I}$ as a one-parameter
family of codifferentials $Q(t)$ of $S(\mathfrak{a})$.

Next we consider the family of maps $\{D_{m_{t}}^{n}\}_{n\in \mathbb{N}}$    %\colon S^{n}(\mathfrak{a})\to \mathfrak{a}\}_{n\in \mathbb{N}}$
defined using the formulae for the higher derived brackets given in \eqref{derivedbrackets}. As explained before,
we can interpret this family of maps as a coderivation of the coalgebra $S(\mathfrak{a}[[t]])$
in the formal setting and
as a one-parameter family of coderivations of the coalgebra $S(\mathfrak{a})$ in the analytical setting. We denote
this coderivation (or family of coderivations respectively) by $M(t)$.

\begin{Lemma}
$M(t)$ satisfies the ordinary differential equation
\begin{equation}
\frac{\dd}{\dd t}Q(t)=M(t)\circ Q(t) - Q(t)\circ M(t).
\label{ODE1}
\end{equation}
\end{Lemma}
%It is clear that the analogous statement is also true in the formal setting.

\begin{proof}
The formula for
$Q(t)$ as a coderivation is %of $S(\mathfrak{a})$ is
\begin{multline*}
Q(t)(x_{1}\otimes \cdots \otimes x_{n})=\\=\sum_{r+s=n}\sum_{\sigma \in (r,s)-\text{shuffles}}{\sign}(\sigma)
D_{E_{t}}^{r}(x_{\sigma(1)}\otimes \cdots \otimes x_{\sigma(r)})
\otimes x_{\sigma(r+1)}\otimes \cdots \otimes x_{\sigma(n)}.
\end{multline*}
As a consequence of \eqref{Xgeneratingphi} we obtain
\begin{multline*}
\frac\dd{\dd t}{Q}(t)(x_{1}\otimes \cdots \otimes x_{n})=\\=\sum_{r+s=n}\sum_{\sigma \in (r,s)-\text{shuffles}}{\sign}(\sigma)
D_{[m_{t},E_{t}]}^{r}(x_{\sigma(1)}\otimes \cdots \otimes x_{\sigma(r)})
\otimes x_{\sigma(r+1)}\otimes \cdots \otimes x_{\sigma(n)}.
\end{multline*}

It is convenient to introduce an auxiliary parameter $\tau$ of degree $1$ and
to consider the $\mathbb{R}[\tau]/\tau^{2}$\ndash modules $\mathfrak{h}[[t]][\tau]/\tau^{2}$ and
$\mathfrak{h}[\tau]/\tau^{2}$. We extend
the graded Lie bracket linearly  by the rule $[\tau x,y]=\tau [x,y]$.
{}From Voronov's result \eqref{J=D} it follows that
\begin{align*}
J_{E_{t}+\tau m_{t}}^{n}=D^{n}_{(E_{t}+\tau m_{t})\circ(E_{t}+\tau m_{t})}=D^{n}_{\tau[m_{t},E_{t}]}=\tau D^{n}_{[m_{t},E_{t}]}.
\end{align*}
Therefore the family of maps $\{(\frac{\partial}{\partial \tau}|_{\tau=0}J^{n}_{E_{t}+\tau m_{t}})\}_{n\in \mathbb{N}}$
corresponds to the coderivation $\dot{Q}(t)$. We claim that $M(t) \circ Q(t) - Q(t) \circ M(t)$ also corresponds to
\begin{align*}
\{(\frac{\partial}{\partial \tau}|_{\tau=0}J^{n}_{E_{t}+\tau m_{t}})\}_{n\in \mathbb{N}},
\end{align*}
which proves the Lemma. 
To verify the claim
it suffices to use the definition  \eqref{Jacobiators} of the Jacobiators,
\begin{multline*}
(J_{E_{t}+\tau m_{t}}^{n})(x_{1}\otimes \cdots\otimes x_{n})=\\=\!\!\!\sum_{r+s=n}\sum_{\sigma \in (r,s)-\text{shuffles}}\!\!\!\!\!\!\!\!\!\!\!\!\sign(\sigma)
D_{E_{t}+\tau m_{t}}^{s+1}(D_{E_{t}+\tau m_{t}}^{r}(x_{\sigma(1)}\otimes 
\cdots \otimes x_{\sigma(r)})\otimes x_{\sigma(r+1)}\otimes \cdots \otimes x_{\sigma(n)}),
\end{multline*}
and to compute
\begin{multline*}
\left.(\frac{\partial}{\partial \tau}\right|_{\tau=0}J_{E_{t}+\tau m_{t}}^{n})(x_{1}\otimes \cdots\otimes x_{n})=\\
\begin{split}
&=\sum_{r+s=n}\sum_{\sigma \in (r,s)-\text{shuffles}}\!\!\!\!\!\!\!\!\!\!\!\!\sign(\sigma)
D_{m_{t}}^{s+1}(D_{E_{t}}^{r}(x_{\sigma(1)}\otimes
\cdots \otimes x_{\sigma(r)})\otimes x_{\sigma(r+1)}\otimes \cdots \otimes x_{\sigma(n)})\\
&-\sum_{r+s=n}\sum_{\sigma \in (r,s)-\text{shuffles}}\!\!\!\!\!\!\!\!\!\!\!\!\sign(\sigma)
D_{E_{t}}^{s+1}(D_{m_{t}}^{r}(x_{\sigma(1)}\otimes 
\cdots \otimes x_{\sigma(r)})\otimes x_{\sigma(r+1)}\otimes \cdots \otimes x_{\sigma(n)}).
\end{split}
\end{multline*}
It is straightforward to see that the first term corresponds to $M(t) \circ Q(t)$ whereas
the second term corresponds to $-Q(t)\circ M(t)$.
\end{proof}

\subsection{Integration to automorphisms}
We now consider the flow of $M(t)$, namely, the solution to
%We will first treat the problem in the formal setting. We search for a formal power series
%$U(t)\colon S(\mathfrak{a}[[t]]) \to S(\mathfrak{a}[[t]])$ that satisfies
\begin{equation}
\label{ODE2}
\begin{aligned}
\frac{\dd}{\dd t}U(t)&=M(t) \circ U(t),\\
U(0)&=\id.
\end{aligned}
\end{equation}
This is equivalent to the following family of equations on the family of maps
$\{U^{n}(t)\}_{n\in\N}$ corresponding to $U(t)$:
\begin{multline}
\label{EQN}
\frac\dd{\dd t}
{U}^{n}(t)(x_{1}\otimes \cdots \otimes x_{n})=\sum_{\sigma \in \Sigma_{n}}\sign(\sigma)\sum_{k\ge 1}\sum_{l_{1}+\cdots+l_{k}=n}\frac{1}{k!l_{1}!\cdots l_{k}!}\\
D_{m_{t}}^{k}\left(U^{l_{1}}(t)(x_{\sigma(1)}\otimes \cdots \otimes x_{\sigma(l_{1})})\otimes \cdots 
\otimes U^{l_{k}}(t)(x_{\sigma(l_{1}+\dots + l_{(k-1)}+1)}\otimes \cdots \otimes x_{\sigma(n)})\right)
\end{multline}
%for all $n \ge 0$, $x_{1},\dots,x_{n} \in \mathfrak{a}[[t]]$. 
%Moreover $U(0)=\id$
together with the initial conditions
$U^{1}(0)=\id$ and
$U^{n}(0)=0$ for $n\ne 1$. 
\begin{Proposition}
The Cauchy problem \eqref{ODE2} has a unique solution. The solution has the property $U^0\equiv0$.
\end{Proposition}
\begin{proof}
That there exists a unique solution for $U(t)$ in the formal setting is seen as follows: first we assume
that we already found (unique) expressions for $U^{m}(t)$, $m <n$. We want to construct $U^{n}(t)$. We expand
it with respect to the formal parameter $t$: $U^{n}(t):=\sum_{r \ge 0} U^{n}_{r}t^{r}$. Condition $U(0)=\id$
determines the term $U^{n}_{0}$ (it is $0$ for $n \ne 0$ and $\id_{\mathfrak{a}}$ for $n=1$). Next suppose
we know $U^{n}_{v}$ for all $v < w$. If we expand equation \eqref{EQN} with respect to the formal parameter
$t$ and consider the term of order $t^{(w-1)}$ we obtain an explicit expression for $U^{n}_{w}$ in terms
of $U^{m}$ for $m<n$ and $U^{n}_{v}$ for $v<w$. So $U^{n}_{w}$ is uniquely determined by these factors. Hence
we can find uniquely determined $U^{n}_{w}$ for all $w \ge 0$ successively and consequently construct $U^{n}$.
We remark that assumption \eqref{formalASS_X} implies $U^{0}(t)=0$.
This completes the proof in the formal setting.

In the analytical setting we first assume that
we have found a family of automorphisms $U(t)\colon S(\mathfrak{a}) \to S(\mathfrak{a})$ integrating the
one-parameter family of coderivations $M(t)$, i.e., solving equation \eqref{ODE2} for $t\in I$.
As before, equation \eqref{ODE2} is equivalent to the family of equations \eqref{EQN}
for all $n \ge 0$, $x_{1},\dots,x_{n} \in \mathfrak{a}$ and $t \in I$.  Moreover $U(t)=\id$
is equivalent to $U^{1}(0)=\id_{\mathfrak{a}}$ and
$U^{n}(0)=0$ for $n > 1$. By assumption \eqref{ASS_X} we can consistently set $U^{0}(t)=0$.

Using uniqueness of solutions of \eqref{ASS_unique} one deduces that a solution to \eqref{ODE2} with $U^{0}(t)=0$
is unique, too:
Suppose we have two solutions satisfying \eqref{ODE2} given by the family of maps $\{U^{n}(t)\}_{n\ge 1}$ and
$\{\tilde{U}^{n}(t)\}_{n\ge 1}$. We consider $\delta U^{n}(t):=U^{n}(t)-\tilde{U}^{n}(t)$. It follows
that $\delta U^{1}(t)$ satisfies \eqref{ASS_unique}, hence $U^{1}(t)=\tilde{U}^{1}(t)$. Now assume
we know that $U^{k}(t)=\tilde{U}^{k}(t)$ for all $k < n$. Equation \eqref{EQN} implies that
$\delta U^{n}(t)$ satisfies \eqref{ASS_unique} too, so $U^{n}(t)=\tilde{U}^{n}(t)$. By induction if follows
that the two solutions coincide.
It remains to prove that such a family of automorphisms $U(t)$ exists for all $t \in I$ under the condition \eqref{ASS_X}.
We inductively define a family of maps 
\[
\{U^{n}(t)\colon S^{n}(\mathfrak{a}) \to \mathfrak{a}\}_{n \ge 1}
\]
that corresponds to an automorphism of $S(\mathfrak{a})$ that satisfies \eqref{ODE2}. (From now on we will suppress the $t$ dependence
of the maps $U^{n}(t)$ and simply write $U^{n}$ instead.)
For $n=1$ we define
\begin{align}
\label{U^{1}}
U^{1}(x):=\Pi_{\mathfrak{a}}\phi_{t}(x).
\end{align}
For $n \ge 1$ we set
\begin{multline}
\label{U}
U^{n}(x_{1}\otimes \cdots \otimes x_{n}):=\sum_{\sigma \in \Sigma_{n}}\sign(\sigma)\sum_{k\ge 1}\sum_{\mu_{1}+\cdots+\mu_{k}=n-1}
\frac{1}{n k! \mu_{1}!\cdots \mu_{k}!}\\
\Pi_{\mathfrak{a}}[[\cdots[\phi_{t}(x_{\sigma(1)}),U^{\mu_{1}}(x_{\sigma(2)}\otimes \cdots \otimes x_{\sigma(\mu_{1}+1)})],\cdots],\\
U^{\mu_{k}}(x_{\sigma(\mu_{1}+\cdots+\mu_{(k-1)}+2)}\otimes \cdots \otimes x_{\sigma(n)})]
\end{multline}
By this formula $U^{n}$ is defined recursively for all $n \ge 1$.

\begin{Lemma}
\label{PropositionX}
The family of maps $\{U^{n}\colon S^{n}(\mathfrak{a})\to \mathfrak{a}\}_{n \ge 1}$ defined by \eqref{U^{1}} and
\eqref{U} satisfies equation \eqref{EQN}
for all $n \ge 1$, $x_{1},\dots,x_{n} \in \mathfrak{a}$ and $t \in I$.
\end{Lemma}
The proof is in the Appendix.
That $U(0)=\id_{S(\fra)}$ can be seen easily: First observe that $U^{1}=\id_{\mathfrak{a}}$ for $t=0$. Moreover,
all $U^{n}$ for $n>1$ vanish at $t=0$ since each term contains the Lie bracket between two elements of $\mathfrak{a}$ which
is an abelian Lie subalgebra. %That equation \eqref{EQN} holds is proven in Appendix A.
\end{proof}
Using equation \eqref{ODE1} one easily deduces that $Z(t):=Q(t)\circ U(t) - U(t)\circ Q(0)$ satisfies
\begin{equation}
\label{ODE3}
\begin{aligned}
\frac{\dd}{\dd t}Z(t)&=M(t) \circ Z(t),\\
Z(0)&=0.
\end{aligned}
\end{equation}
In the formal setting one immediately proves that $Z(t)=0$ is the unique solution to \eqref{ODE3} (under assumption \eqref{formalASS_X}).
In the analytical setting one first computes $Z^{0}=\Pi_{\mathfrak{a}}\phi_{t}P-\Pi_{\mathfrak{a}}\phi_{t}\Pi_{\mathfrak{a}}P$
(recall that our Maurer--Cartan element is $E=\Hat E + [P,\cdot]$)
%$P$ denotes the element such that $E=E_{1}+[P,\cdot]$) 
which vanishes because of \eqref{globalASS_X}. Now one can apply the same arguments as in the proof of uniqueness for $U(t)$
and one obtains that $Z(t)=0$.

By definition of $Z(t)$, $Z(t)=0$ is equivalent to
\begin{align}
Q(t)\circ U(t) = U(t) \circ Q(0)
\end{align}
which means that $U(t)$ defines an $L_\infty$\ndash isomorphism.
This completes the proof of Theorems~\ref{formal} and~\ref{MAIN}.

\section{Applications}\label{s:app}% in Poisson Geometry}

We describe an application of Theorem~\ref{MAIN} in the framework of Poisson geometry. Out of it applications in symplectic and Dirac geometry follow.

Let $M$ be a smooth finite-dimensional manifold.
As noticed
by Oh and Park in \cite{OhPark}, if $M$ is a Poisson manifold,
the space of sections of the exterior algebra of the normal bundle of a submanifold of a certain class
(namely, a coisotropic submanifold) carries the structure of a flat $L_{\infty}$\ndash algebra. The
same structure was found in \cite{CattaneoFelder}
as the semi-classical limit of a certain topological quantum field theory called the Poisson Sigma model; 
the $L_{\infty}$\ndash algebra structure was derived
not only for coisotropic submanifolds but for every submanifold of $M$ (coisotropic submanifolds are special in so far
as they are exactly those whose associated $L_\infty$\ndash algebras are flat).
We now briefly recall the construction in \cite{CattaneoFelder}, which
makes use of graded manifolds and Voronov's higher derived brackets. 
%We now recall the version of \cite{C} which does not use graded manifolds.

\subsection{Submanifolds and V-algebras}
Given a smooth manifold $M$,
the space of multivector fields
$\mathcal{V}(M)[1]:=\Gamma(M,\Lambda TM)[1]$ carries the structure of a graded Lie algebra where
the graded Lie bracket is given by the Schouten--Nijenhuis bracket which we denote by $[\cdot,\cdot]$; see Example~\ref{exaDerA}.

Let $S$ be a submanifold. Its normal bundle $NS$ is by definition the quotient of the restriction $T_SM$ of $TM$ to $S$ by $TS$.
Set $A:=\Gamma(S,\Lambda NS)$ and $\fra:=A[1]$ as in Example~\ref{exaDerA}.
By restricting a multivector field to $S$ and then projecting it to its normal
components, we get a projection $\Pi_{M;\fra}\colon\mathcal{V}(M)[1]\to\fra$. 
%for which sections do not exist in general. They exist however on the quotient 
%of $\mathcal{V}(M)$ consisting of germs of multivector fields along $S$.
%Poisson bivector fields on $M$ corresponds exactly to Maurer--Cartan elements of $(\mathcal{V}(M)[1],[\cdot,\cdot])$.
%Let $S$ be the submanifold under consideration and
%Let $\sim_{S}$ be the following equivalence relation on $\mathcal{V}(M)$:
%$X \sim_{S} Y$ if and only if there exists an open neighbourhood $U$ of $S$ such that $X|_{U}=Y|_{U}$. We set
%$\mathcal{V}(M,S):=\mathcal{V}(M)/\sim_{S}$, the Gerstenhaber algebra of germs of multivector fields.
%The space
%$\frh_{M,S}:=\mathcal{V}(M,S)[1]$ inherits both the structure of a graded Lie algebra and a projection $\Pi_{M,S;\fra}$
%onto $\fra$. %As sections now exist, $\frh_{M,S}$ can be made into a V\ndash algebra.
%As we will shortly see, it also has the structure of a V\ndash algebra though not in a canonical way.
%upon choosing an embedding $\sigma\colon NS\hookrightarrow M$ with $\sigma(S)=S$.

Denote the vanishing ideal of $S$ by
$I(S):=\{f\in\mathcal{C}^{\infty}(M)| \; f|_{C}=0\}$.
The inclusions $i_{n m}: I^{m}(S) \hookrightarrow
I^{n}(S)$ for $m \ge n$ equip the collection $\mathcal{V}(M)/I^{n}(S)\mathcal{V}(M)$ with the structure of a
projective system and we define the Gerstenhaber algebra of multivector fields on a formal neighbourhood of $S$ in $M$ by
\begin{align*}
\mathcal{V}(M,S):=\lim_{\leftarrow}\mathcal{V}(M)/I^{n}(S)\mathcal{V}(M).
\end{align*}
The space $\frh_{M,S}:=\mathcal{V}(M,S)[1]$ inherits both the structure of a graded Lie algebra and a
projection $\Pi_{M,S;\fra}$ onto $\fra$.
As we will shortly see, it also has the structure of a V\ndash algebra though not in a canonical way.

Thus, a Maurer--Cartan element of $\frh_{M,S}$ induces an $L_\infty$\ndash structure on $\fra$.
Observe that the class $[\pi]$ in $\frh_{M,S}$ of a bivector field $\pi$ on $M$ is a Maurer--Cartan element if and only if the restrictions to $S$ of
$[\pi,\pi]$ and all its derivatives vanish. In this case, we say that $\pi$ is Poisson in a formal neighbourhood of $S$.
Moreover, $\Pi_{M,S;\fra}[\pi]$ vanishes if and only if $\pi_x(\alpha,\beta)=0$ $\forall x\in S$, $\forall\alpha,\beta\in N^*_xS$. In this
case $S$ is called a coisotropic submanifold.

We now explain how to induce a V\ndash structure on $\frh_{M,S}$ using a choice of embedding $\sigma\colon NS\hookrightarrow M$
with $\sigma|_{S}=\id_{S}$.
%However, to apply Theorem~\ref{MAIN}, we need a fix V\ndash algebra structure. To this aim, we have to resort to the description given in \cite{CattaneoFelder}.
%Let $\frp:=\Ker\Pi_\fra$ and $\frh:=\fra\oplus\frp$. We have an isomorphism 
Regard $A$ as a graded commutative algebra and set $\frh:=\Hat S_A(\Der(A)[-1])[1]$ with the V\ndash algebra structure
of Example~\ref{exaDerA}. We now claim that $\frh$ is isomorphic, though noncanonically, to $\frh_{M,S}$.
To do this, we observe that
$A$ is the algebra of functions on the graded manifold $N^*[1]S$. %where the conormal bundle $N^*S$ is by definition the kernel of the projection $T^*_SM\to T^*S$.
So $\frh[-1]$ is the formally completed
Gerstenhaber algebra of multivector fields on $N^*[1]S$.
By  the Legendre mapping theorem \cite{Roytenberg}, this is canonically isomorphic to the formally completed Gerstenhaber algebra of multivector
fields on the graded manifold $N[0]S$ which is the same as the Gerstenhaber algebra $\calV(NS,S)$ of multivector fields on a formal neighbourhood
of $S$ in $NS$.
Finally, the choice of embedding $\sigma$ yields an isomorphism between 
$\calV(NS,S)$ and $\calV(M,S)$, and so an isomorphism $\Hat\sigma\colon\frh\to\frh_{M,S}$.

Two different choices of embeddings yield an automorphism of the graded Lie algebra $\frh$. We will see in the next subsection that the assumption
of Theorem~\ref{MAIN} are respected, so the effect of a change of embedding may be understood easily now.

\begin{Remark}
A simpler construction, avoiding graded manifolds, is that of \cite{C}. It starts with the observation that an embedding $\sigma$ 
yields a section $\Tilde\sigma\colon\fra\to\frh_{M,S}$ with the property
% However, this construction is not suitable for discussing the changes of embeddings.
% Namely,
% a section $\Tilde\sigma\colon\fra\to\frh_{M,S}$ can be induced from the choice of embedding. % $\sigma\colon NS\hookrightarrow M$ with $\sigma(S)=S$
% Namely, let $U:=\sigma(NS)$ and $p_U\colon U\to S$ the induced projection. If $f$ is a function on $S$, we define $\Tilde\sigma f:=[p^*f]$,
% where $[\ ]$ denotes the class in $\mathcal{V}(M,S)$. Next observe that $\sigma$ induces a decomposition $T_SU=TS\oplus NS$.
% So a section  $X$ of $NS$ may be regarded as a section of $T_SU$ and lifted to a vertical vector field $\Tilde X$ on $U$ which is constant along the fibers.
% We then set $\Tilde\sigma(X)=[\Tilde X]$. Finally, we extend $\Tilde\sigma$ to the whole of $\fra$ as an algebra homomorphism. 
%It is not difficult to check 
that $\Tilde\sigma(\fra)$ is an abelian subalgebra.
Let $\frp:=\Ker\Pi_{M,S;\fra}$ and $\iota_\frp$ its inclusion map into $\frh_{M,S}$.
We then have the isomorphism $\Tilde\sigma\oplus\iota_\frp\colon\fra\oplus\frp\to\frh_{M,S}$. This induces a V\ndash algebra structure on 
$\fra\oplus\frp$. Notice however that the Lie bracket on $\fra\oplus\frp$ depends on the choice of embedding.
Hence this simpler construction, while perfectly fine for inducing $L_\infty$\ndash structures on $\fra$, is not suitable for the application of Theorem~\ref{MAIN}
and so for discussing the effect of  a change of embedding.
% Observe that, given a Maurer--Cartan element in $\frh_{M,S}$,
% the induced $L_\infty$\ndash structure on $\fra$ depends on the choice of embedding of $NS$ as both the Maurer--Cartan element and
% the Lie bracket on $\fra\oplus\frp$ depend on this choice. 
\end{Remark}

\begin{Remark}
As already remarked, the induced $L_\infty$\ndash structure is flat if and only if
$S$ is a coisotropic submanifold. In this case, one can show \cite{OhPark,CattaneoFelder,C} that the unary operation does not depend on the choice
of embedding and is the Lie algebroid differential associated to the conormal bundle of $S$ as a Lie subalgebroid of the cotangent bundle of $M$.
\end{Remark}

\subsection{Uniqueness of the induced $L_\infty$-structure}
It is well-known from differential topology (see \cite{Hirsch} for instance) that any two
tubular neighbourhoods of $S$ in $M$ are isotopic. For our purposes this can be expressed as follows: For any two embeddings $\sigma_0$ and $\sigma_1$
of $NS$ into $M$,
there is a family $V_{t}$, $t \in I=[0,1]$,
of open neighbourhoods of $S$ in $M$,
a family of diffeomorphisms $\psi_{t}\colon V_0 \to V_t$ and a family of embeddings $\sigma_t\colon NS\to M$,
such that $\psi_{0}=\id_{V_0}$, $\psi_{t}|_{S}=\id_{S}$,
and $\psi_t\circ\sigma_0=\sigma_t$ in an open neighbourhood of $S$. %We define $\sigma_t:=\psi_t^{-1}\circ\sigma_0$.
The pushforward $\psi_{t*}$ of multivector fields defines an automorphism of $\frh_{M,S}$ which we denote by $\Hat\psi_t$.
Denoting by $\Hat\sigma_t$ the isomorphism $\frh\to\frh_{M,S}$ induced by $\sigma_t$, we then get
$\Hat\psi_t\circ\Hat\sigma_0=\Hat\sigma_t$.
Let $\phi_t:=\Hat\sigma_t^{-1}\circ\Hat\sigma_0=\Hat\sigma_0^{-1}\circ\Hat\psi_t^{-1}\circ\Hat\sigma_0$.
Let $Z_t:=-\frac\dd{\dd t}\psi_t$ as a vector field in an open neighbourhood of $S$ and $\Hat Z_t$ its class in $\frh_{M,S}$.
Then equation \eqref{Xgeneratingphi} is satisfied with $m_t=[\Hat\sigma_0^{-1}(\Hat Z_t),\cdot]$. Observe that $Z_t|_S$ is tangent to $S$.
Using the explicit formula for the Legendre mapping, it is easy to verify that this implies condition \eqref{ASS_X}.
% %and $i_{0}\circ \psi_{1}=i_{1}$ (restricted to $U$).
% \\
% From this we obtain a family of automorphisms $\hat{\psi}^{*}_{t}$ of the algebra $A$ and a family of 
% automorphisms $\phi_{t}$ of $\mathcal{A}$ generated by an inner derivation $[X_{t},\cdot]$ with $X_{t} \in \mathcal{A}^{0}$.
% $X_{t}$
% can be thought of the jet of vector fields in transversal directions generating $\hat{\psi}^{*}_{t}$.
% \\
% With help of the two two embeddings $i_{0}$ and $i_{1}$ one constructs two isomorphisms
% $\varphi_{0}$, $\varphi_{1}: \mathcal{V}(M,S) \to \mathcal{A}$. The fact that
% $i_{0}\circ \psi_{1}=i_{1}$ implies 
% \begin{align}
% \label{relation}
% \phi_{1}\circ \varphi_{0}=\varphi_{1}.
% \end{align}
% Define $P:=\varphi_{0}(\hat{\Pi})$ and $P_{t}:=\phi_{t}(P)$ which is a family of Maurer--Cartan elements
% of $\mathcal{A}$. By equation \eqref{relation} one has $\varphi_{1}(\hat{\Pi})=P_{1}$, so
% we obtain a one-parameter family of $L_\infty$\ndash algebra structures on $\Gamma(S,\bigwedge NS)$
% that connects the $L_\infty$\ndash algebra structures constructed with help of the embeddings $i_{0}$ and $i_{1}$
% respectively.
% \\
Finally, uniqueness of solutions
of equation \eqref{ASS_unique} follows from the
uniqueness of flows generated by vector fields on graded manifolds (in view of the canonical isomorphism between $\frh$ and $\calV(NS,S)$
this is in this case just the uniqueness of flows generated by vector fields on $NS$).
So all assumptions of Theorem~\ref{MAIN} hold and one concludes:

\begin{Theorem}
\label{intrinsic}
The $L_{\infty}$\ndash algebra structures constructed on $\fra$ with the help of two different
embeddings of $NS$ into $M$ as tubular neighbourhoods of $S$ are $L_{\infty}$\ndash isomorphic.
\end{Theorem}

\begin{Remark}\label{Remark1}
In case one changes the tubular neighbourhood by acting on $NS$ via
a vector bundle automorphism, there is a simpler proof by
applying the construction in Example~\ref{Example1}: in fact the vector bundle automorphism induces
an automorphism of $A:=\Gamma(S,\Lambda NS)$, and the natural extension to an automorphism
of $\frh:=\Hat S_A(\Der(A)[-1])[1]$ also relates the two associated Maurer--Cartan elements.
Consequently the induced $L_{\infty}$\ndash algebras on $\fra:=A[1]$ are $L_{\infty}$\ndash isomorphic,
and the $L_\infty$\ndash isomorphism is linear.
\end{Remark}

Theorem~\ref{intrinsic} immediately implies
the following

\begin{Corollary}
\label{extrinsic}
Let $(M_{1},\pi_{1})$ and $(M_{2},\pi_{2})$ be two Poisson manifolds and $S_{1}$, $S_{2}$ submanifolds
of $M_{1}$ and $M_{2}$ respectively. Assume
$\psi\colon  M_{1} \to M_{2}$ is a Poisson diffeomorphism that maps $S_{1}$ to $S_{2}$. Then the
isomorphism classes of the two $L_\infty$\ndash algebras associated to $S_{1}$ and $S_{2}$ coincide.
\end{Corollary}

\begin{proof}
Fix an embedding of $NS_{1}$ into $M_{1}$. The diffeomorphism $\psi$ induces a bundle isomorphism between
$NS_{1}$ and $NS_{2}$ and using this identification we obtain an embedding of $NS_{2}$ into $M_{2}$.
Hence $\psi$ allows us to identify the two V\ndash algebras associated to $S_{1}$ and $S_{2}$. Moreover,
the Maurer--Cartan elements associated to $\pi_{1}$ and to $\pi_{2}$ also get identified via $\psi$. 
So the two induced $L_{\infty}$\ndash algebras are $L_{\infty}$\ndash isomorphic. By Theorem~\ref{intrinsic}
other choices of embeddings of $NS_{1}$ and $NS_{2}$ into $M_{1}$ and $M_{2}$, respectively, will not affect the isomorphism
classes of the two $L_{\infty}$\ndash algebras.
\end{proof}

\subsection{Presymplectic manifolds}
Let $S$ be a finite-dimensional smooth manifold.
A two-form $\omega$ on $S$ may be regarded as a bundle map $\omega^\sharp\colon TS\to T^*S$
by $\omega^\sharp_x(v):=\omega_x(v,\ )$. If $\omega$ is closed and $\omega^\sharp$ has constant rank,
$S$ is called a presymplectic manifold. If the rank is maximal (i.e., $\omega^\sharp$ is bijective), 
then $S$ is called a symplectic manifold.
A symplectic manifold is also a Poisson manifold with Poisson bivector field obtained by inverting the symplectic two-form.
A coisotropic submanifold in a symplectic manifold gets the structure of a presymplectic submanifold by restricting the
symplectic form. 

Let $(S,\omega)$ be a presymplectic manifold. Then $\calF_\omega:=\Ker\omega^\sharp$ is an integrable distribution.
Thus, the de~Rham differential descends to the quotient 
\[
\Omega_{\calF_\omega}:=\Omega(S)/\{\alpha\in\Omega(S):i_X\alpha=0\ \forall X\in\Gamma(S,\calF_\omega)\}
\]
called the foliated de~Rham complex. Also observe that
$\Omega_{\calF_\omega}=\Gamma(S,\Lambda\calF^*_\omega)$.

\begin{Corollary}[Oh--Park]\label{cor:pres}
The foliated de~Rham complex $\Omega_{\calF_\omega}$ of a presymplectic manifold $(S,\omega)$ carries a flat $L_\infty$\ndash structure,
unique up to $L_\infty$\ndash automorphisms, with first operation the de~Rham differential.
\end{Corollary}

See \cite{OhPark} for a different proof.
%This Corollary is already contained in \cite{OhPark}, where however the uniqueness part of the proof is incomplete.

\begin{proof}
By a theorem of Gotay \cite{G}, every presymplectic manifold $(S,\omega)$ may be embedded into some symplectic manifold
$(M,\Omega)$ as a coisotropic submanifold with $\omega=\iota^*\Omega$, where $\iota\colon S\to M$ is the embedding.
Moreover, $\Omega^\sharp$ establishes an isomorphism of $\calF_\omega$ with $N^*S$. So the construction in the first part of this Section
yields the desired flat $L_\infty$\ndash structure.

Gotay also proves that this coisotropic embedding is unique up to neighbourhood equivalence: namely,
for every two coisotropic embeddings of $S$, there exist symplectomorphic neighbourhoods of $S$.
Applying Corollary~\ref{extrinsic}, we get uniqueness.
%$\iota_i\colon S\to M_i$, $i=,12$ with $\omega=\iota_1^*\Omega_1=\iota_2^*\Omega_2$,
%there exist symplectomorphic neighbourhoods of $S$ 
\end{proof}

\subsection{Regular Dirac structures}
Let $S$ be a smooth manifold. Sections of $TS\oplus T^*S$ may be endowed with the Courant bracket \cite{Cour} which is
the skew--symmetrization of the Dorfman bracket \cite{Dorf} given by
\[
[X_1\oplus\xi_1,X_2\oplus\xi_2]=
[X_1,X_2]\oplus\left(
L_{X_1}\xi_2 - i_{X_2}d\xi_1
\right)
\]
and with the symmetric nondegenerate pairing 
$\langle X_1\oplus\xi_1,X_2\oplus\xi_2\rangle = i_{X_1}\xi_2 + i_{X_2}\xi_1$.
A subbundle $L$ of $TS\oplus T^*S$ is called a Dirac structure if it is maximally isotropic with respect to the pairing
and sections of $L$ are closed under the Courant bracket.
Examples of Dirac structures are graphs of Poisson bivector fields.

A Dirac structure $(S,L)$ is called regular if $\calF_L:=L\cap TS$ has constant rank.
Examples of regular Dirac structures are graphs of presymplectic forms.
Coisotropic submanifolds of a Poisson manifold with regular characteristic distribution
get an induced regular Dirac structure.
Since $\calF_L$ is an integrable distribution, one can define the foliated de~Rham complex $\Omega_{\calF_L}$. 
%as in the case of presymplectic manifold. 
We then have the following generalization of Corollary~\ref{cor:pres}:

\begin{Corollary}
The foliated de~Rham complex $\Omega_{\calF_L}$ of a regular Dirac manifold $(S,L)$ carries a flat $L_\infty$\ndash structure,
unique up to $L_\infty$\ndash automorphisms, with first operation the de~Rham differential.
\end{Corollary}

Notice that the existence part is already contained in \cite{CZ}.

\begin{proof}
It is shown in \cite{CZ} that, canonically up to neighbourhood equivalences,
the total space of $\calF_L^*$ can be given a Poisson structure such that the zero section is coisotropic
with induced Dirac structure equal to $L$. In particular the Poisson structure establishes an isomorphism $N^*S\to\calF_L$.
\end{proof}

%\begin{appendix}
\appendix

\section*{Appendix. Proof of Lemma~\ref{PropositionX}}

We prove that the family of maps $\{U^{n}\colon S^{n}(\mathfrak{a})\to \mathfrak{a}\}_{n \ge 1}$
defined by equations \eqref{U^{1}} and \eqref{U}
satisfies the family of relations given
by equation \eqref{EQN} --- again we suppress the $t$ dependence of $U^{n}(t)$.
The proof we give works inductively: It is easy to check
that $U^{1}(a_{1}):=\Pi_{\mathfrak{a}}\phi_{t}(a_{1})$ satisfies
$\dot{U}^{1}=\Pi_{\mathfrak{a}}m_{t}\circ U^{1}$, which
is equation \eqref{EQN} for $n=1$.

Suppose we verified that equation \eqref{EQN} holds for all $U^{k}$, $k < n$. We show that
this implies that equation \eqref{EQN} is satisfied for $n$, too. The definition of $U^{n}$ by equation \eqref{U}
implies
\begin{multline*}
\dot{U}^{n}(a_{1}\otimes \cdots \otimes a_{n})=\sum_{\sigma \in \Sigma_{n}}\sign(\sigma)\sum_{k\ge 1}\sum_{\mu_{1}+\cdots+\mu_{k}=n-1}
\frac{1}{n k! \mu_{1}!\cdots \mu_{k}!}\\
\left(\Pi_{\mathfrak{a}}[[\cdots[[m_{t}\phi_{t},U^{\mu_{1}}],U^{\mu_{2}}],\cdots],U^{\mu_{k}}]+
k\Pi_{\mathfrak{a}}[[[\cdots[\phi_{t},U^{\mu_{1}}],\cdots ],U^{\mu_{(k-1)}}],\dot{U}^{\mu_{k}}]\right),
\end{multline*}
where we suppressed the arguments $(a_{\sigma(1)},\dots,a_{\sigma(n)})$. The first term
comes from deriving $\phi_{t}$, the second one from
deriving one of the factors $U^{k}$ with $k <n$ in the formula for $U^{n}$.
We denote the two terms by $A^{n}$ and $B^{n}$ respectively. $A^{n}$ contains
terms of the form $[[[m_{t}\phi_{t},U^{\mu_{1}}],U^{\mu_{2}}],\dots]$ where we can first use that $m_{t}$ is a derivation and then
successively apply the graded Jacobi identity
(see Definition~\ref{DGLA}) and obtain
\begin{multline*}
A^{n}=\sum_{\sigma \in \Sigma_{n}}\sign(\sigma)\sum_{k\ge 1}\sum_{r+s=k}
\sum_{\alpha_{1}+\cdots+\alpha_{r}+ \atop +\beta_{1}+\cdots +\beta_{s}=n-1}
\frac{1}{n r! s! \alpha_{1}!\cdots \alpha_{r}! \beta_{1}! \cdots \beta_{s}!} \\
\Pi_{\mathfrak{a}}[([[\cdots[m_{t}U^{\alpha_{1}},U^{\alpha_{2}}],\cdots ],U^{\alpha_{r}}]),([[\cdots[\phi_{t},U^{\beta_{1}}],\cdots],U^{\beta_{s}}])].
\end{multline*}
Next we apply equation \eqref{projection} which leads to
\begin{multline*}
A^{n}=\bigg(\sum_{\sigma \in \Sigma_{n}}\sign(\sigma)\sum_{k\ge 1}\sum_{r+s=k}
\sum_{\alpha_{1}+\cdots+\alpha_{r}+ \atop +\beta_{1}+\cdots +\beta_{s}=n-1}
\frac{1}{n r! s! \alpha_{1}!\cdots \alpha_{r}! \beta_{1}! \cdots \beta_{s}!} \\
D_{m_{t}}^{r+1}(U^{\alpha_{1}}\otimes \cdots \otimes
U^{\alpha_{r}}\otimes \Pi_{\mathfrak{a}}[[\cdots[\phi_{t},U^{\beta_{1}}],\cdots ],U^{\beta_{s}}]\bigg)-\\
-\bigg(\sum_{\sigma \in \Sigma_{n}}\sign(\sigma)\sum_{k\ge 1}\sum_{r+s=k}
\sum_{\alpha_{1}+\cdots+\alpha_{r}+ \atop +\beta_{1}+\cdots +\beta_{s}=n-1}
\frac{1}{n r! s! \alpha_{1}!\cdots \alpha_{r}! \beta_{1}! \cdots \beta_{s}!} \\
\Pi_{\mathfrak{a}}[([[\cdots[\phi_{t},U^{\alpha_{1}}],\cdots ],U^{\alpha_{r}}]),\Pi_{\mathfrak{a}}([[\cdots[m_{t}U^{\beta_{1}},U^{\beta_{2}}],\cdots],U^{\beta_{s}}])]
\bigg).
\end{multline*}
We claim that the following two identities hold: the first is
\begin{multline*}
\sum_{\sigma \in \Sigma_{n}}\sign(\sigma)\sum_{k\ge 1}\sum_{\mu_{1}+\cdots+\mu_{k}=n-1}
\frac{1}{n (k-1)! \mu_{1}!\cdots \mu_{k}!}\\
\Pi_{\mathfrak{a}}[[[\cdots[\phi_{t},U^{\mu_{1}}],\cdots ],U^{\mu_{(k-1)}}],\dot{U}^{\mu_{k}}])=\\
=\sum_{\sigma \in \Sigma_{n}}\sign(\sigma)\sum_{k\ge 1}\sum_{r+s=k}
\sum_{\alpha_{1}+\cdots+\alpha_{r}+ \atop + \beta_{1}+\cdots +\beta_{s}=n-1}
\frac{1}{n r! s! \alpha_{1}!\cdots \alpha_{r}! \beta_{1}! \cdots \beta_{s}!} \\
\Pi_{\mathfrak{a}}[([[\cdots[\phi_{t},U^{\alpha_{1}}],\cdots ],U^{\alpha_{r}}]),\Pi_{\mathfrak{a}}([[\cdots[m_{t}U^{\beta_{1}},U^{\beta_{2}}],\cdots],U^{\beta_{s}}])],
\end{multline*}
which means that $B^{n}$ cancels with the second term in the expression for $A^{n}$ given above; the second is
\begin{multline*}
\bigg(\sum_{\sigma \in \Sigma_{n}}\sign(\sigma)\sum_{k\ge 1}\sum_{r+s=k}
\sum_{\alpha_{1}+\cdots+\alpha_{r}+ \atop +\beta_{1}+\cdots +\beta_{s}=n-1}
\frac{1}{n r! s! \alpha_{1}!\cdots \alpha_{r}! \beta_{1}! \cdots \beta_{s}!} \\
\hspace{3cm} \left(D_{m_{t}}^{r+1}(U^{\alpha_{1}}\otimes \cdots \otimes
U^{\alpha_{r}}\otimes \Pi_{\mathfrak{a}}[[\cdots[\phi_{t},U^{\beta_{1}}],\cdots ],U^{\beta_{s}}])\right)\bigg)\\
=\sum_{\sigma \in \Sigma_{n}}\sign(\sigma)\sum_{k\ge 1}\sum_{l_{1}+\cdots+l_{k}=n}
\frac{1}{k! l_{1}!\cdots l_{k}!}
D_{m_{t}}^{k}(U^{l_{1}}\otimes \cdots \otimes U^{l_{k}})
\end{multline*}
which means that the first term in the expression for $A^{n}$
is equal to the expression from equation \eqref{EQN} which we would like to obtain.

The first identity is straightforward to check: By the induction hypothesis, equation \eqref{EQN}
is satisfied for $k <n$, so we can plug in the expression for $\dot{U}^{\mu_{k}}$ on
the left-hand side of the identity. This immediately leads to the expression on the right-hand side.
To prove the second identity, we first use the recursive definition of $U^{n}$ (see formula
\eqref{U}) on the left-hand side of the identity to arrange
the terms of the form $\Pi_{\mathfrak{a}}[[\cdots[\phi_{t},U^{\beta_{1}}],\cdots ],U^{\beta_{s}}]$
into some $U^{\beta}$. We arrive at
\begin{align}
\label{2ndsum}
\sum_{\sigma \in \Sigma_{n}}\sign(\sigma)\sum_{r\ge 1}\sum_{\alpha_{1}+\cdots+\alpha_{r}=n}
\frac{1}{n (r-1)! (\alpha_{1}-1)!\alpha_{2}!\cdots \alpha_{r}!}
D_{m_{t}}^{r}(U^{\alpha_{1}}\otimes \cdots \otimes U^{\alpha_{r}}).
\end{align}
It remains to prove that this map is equal to
\begin{align}
\label{1stsum}
\sum_{\sigma \in \Sigma_{n}}\sign(\sigma)\sum_{k\ge 1}\sum_{l_{1}+\cdots+l_{k}=n}
\frac{1}{k! l_{1}!\cdots l_{k}!}
D_{m_{t}}^{k}(U^{l_{1}}\otimes \cdots \otimes U^{l_{k}}).
\end{align}
We give the construction of a third map for which it is easy to show that it is equal to both map \eqref{2ndsum} and map \eqref{1stsum}.
Assume one is given $n$ distinguishable objects and $r$ `boxes' where
there are $w_{j}$ boxes that can contain exactly $l_{j}$ of the objects, $1\le j \le k$
($0 < l_{1} < \dots < l_{k}$ and $w_{1}+\dots + w_{k}=r$). We label this situation by $(r|(l_{1},w_{1}),\dots,(l_{k},w_{k}))$.
We assume that boxes that contain the same number of objects are indistinguishable.
The number of different ways to put the $n$ objects into these boxes is given by
\begin{align*}
\frac{n!}{w_{1}!\cdots w_{k}!(l_{1}!)^{w_{1}}\cdots (l_{k}!)^{w_{k}}}.
\end{align*}
Consider
\begin{align*}
\sum_{\sigma \in \Sigma_{n}}\sign(\sigma)\sum_{r\ge 1}\sum_{(r|(l_{1},w_{1}),\dots,(l_{k},w_{k}))}
|\text{ways to put $n$ objects into these boxes}|\\
\hspace{3cm}\left(D_{m_{t}}^{r}(U^{l_{1}}\otimes \cdots \otimes U^{l_{1}}\otimes \cdots \otimes U^{l_{k}}\otimes \cdots \otimes U^{l_{k}})\right).
\end{align*}
It is straightforward to check that this map is equal to map \eqref{2ndsum} and to map \eqref{1stsum}.

%\end{appendix}

%\begin{thebibliography}{AAAAA}

%\end{thebibliography}

\end{document}